\documentstyle[12pt,amssymb,amsmath]{article}
\newcounter{ar}

\newcounter{bk}

\newcommand{\mysection}[1]{\section{\large\bf #1}}

\newtheorem{defn}[subsection]{Definition}
\newtheorem{prop}[subsection]{Proposition}
\newtheorem{th}[subsection]{Theorem}

\newtheorem{deflem}[subsection]{Definition-Lemma}
\newtheorem{corll}[subsection]{Corollary}
\newtheorem{conj}[subsection]{Conjecture}

\newenvironment{ssect}[1]{\smallskip\noindent%
\refstepcounter{subsection}%
{\bf \thesubsection~#1}\hspace{-1mm}}
{\smallskip}

\newenvironment{rem}{\smallskip\noindent%
\refstepcounter{subsection}%
{\bf \thesubsection}~~{\sc Remark.}\hspace{-1mm}}
{\smallskip}

\newenvironment{ex}{\smallskip\noindent%
\refstepcounter{subsection}%
{\bf \thesubsection}~~{\sc Example.}}
{\smallskip}

\newenvironment{proof}[1]{\noindent {\em Proof#1.}}
{~$\blacksquare$\smallskip}

\newenvironment{ackn}{\medskip \noindent \small
{\sl Acknowledgments.}}{\bigskip}

\newenvironment{smallbibl}[1]
{\small
\begin{thebibliography}{#1}}
{\end{thebibliography}}

\newcommand{\res}{{\rm res}}

\newcommand{\dvsr}{{\rm div}}
\newcommand{\lk}{\ell k_{polar}}
\newcommand{\C}{{\Bbb C}}
\newcommand{\R}{{\Bbb R}}
\newcommand{\Q}{{\Bbb Q}}
\renewcommand{\P}{{\Bbb P}}

\newcommand{\Cdot}{\,\mbox{{\large$\cdot$}}\,}

\newcommand{\tpi}{2\pi i}

\newcommand{\myfigure}{

\unitlength0.75mm
\begin{picture}(130,60)(-30,-6)

\thicklines
\put(40,0){\line(1,0){70}}
\put(75,0){\circle*{1.5}}
\put(75,10){\line(0,1){37}}

\qbezier(110,20)(-30.3,30)(110,40)
\put(75,37){\circle*{1.5}}
\put(75,23){\circle*{1.5}}

\put(100,8.2){\line(0,1){2.5}}
\put(100,12.1){\line(0,1){2.6}}
\put(100,16.1){\line(0,1){2.6}}
\put(100,20.1){\line(0,1){2.6}}
\put(100,24.1){\line(0,1){2.6}}
\put(100,28.1){\line(0,1){2.6}}
\put(100,32.1){\line(0,1){2.6}}
\put(100,36.1){\line(0,1){2.6}}
\put(100,40.1){\line(0,1){2.6}}
\put(100,44.1){\line(0,1){2}}

\thinlines
\put( 40,10){\line(0,1){37}}
\put(110,10){\line(0,1){37}}

\qbezier(40,47)(65,50)( 75,47)
\qbezier(75,47)(85,44)(110,47)

\qbezier(40,10)(65,15)( 75,10)
\qbezier(75,10)(85, 5)(110,10)

\put(34,39){\makebox(0,0)[b]{\small$A$}}
\put(39,-.5){\makebox(0,0)[br]{\small$f(A)$}}
\put(75,-3){\makebox(0,0)[t]{\footnotesize$Q$}}
\put(52,36){\makebox(0,0)[b]{\footnotesize$N$}}
\put(76,12){\makebox(0,0)[bl]{\footnotesize$C$}}
\put(76,24){\makebox(0,0)[bl]{\footnotesize$P_1$}}
\put(76,39){\makebox(0,0)[bl]{\footnotesize$P_2$}}
\put(101,11){\makebox(0,0)[bl]{\footnotesize$D$}}

\put(34,34){\vector(0,-1){27}}

\end{picture}
}

\begin{document}


\title{
\vspace{-2.3cm}\begin{flushright}
{\normalsize\bf ITEP-TH-46/00}\\
\end{flushright}
\vspace{1.3cm}
{\Large\sc Polar Homology}
}

\author{Boris Khesin\thanks{Department of Mathematics,
University of Toronto,
Toronto, ON M5S 3G3, Canada;
e-mail: {\tt khesin@math.toronto.edu}}~
and Alexei Rosly\thanks{Institute of Theoretical and
Experimental Physics, B.Cheremushkinskaya 25, 117259 Moscow,
\mbox{Russia}; e-mail: {\tt rosly@heron.itep.ru}}}

\date{August 12, 2000}

\maketitle

\begin{abstract}
For complex projective manifolds
we introduce polar homology groups,
which are
holomorphic analogues of the homology groups in topology.
The polar $k$-chains are subvarieties of complex
dimension $k$ with meromorphic forms on them,
while the boundary operator is defined by
taking the polar divisor and the Poincar\'e
residue on it.
\end{abstract}


\mysection{Introduction} \label{Int}

In this paper we  introduce certain homology
groups defined for complex projective manifolds that can be regarded as a
complex version of singular   homology groups in topology. The idea
of such a geometric analogue of topological homology comes from
thinking of the Dolbeault (or $\bar\partial$) complex of $(0,k)$-forms on a
complex manifold as an obvious analogue of the de Rham complex of $k$-forms
on a smooth manifold. This poses an immediate question: ``What is the
analogue of the chain complex relevant to the context of complex
manifolds?", which we address in detail below.

It should be  mentioned that the
correspondence between de Rham and Dolbeault complexes, or
$d\leftrightarrow\bar\partial$, has the following natural extension.
\vspace{-3mm}
$$
\begin{array}{rcl}
                     d & \leftrightarrow & \bar\partial \\
{\rm de~Rham~complex} & \leftrightarrow & {\rm Dolbeault~complex} \\
{\rm smooth~functions~or~sections} & \leftrightarrow &
                              {\rm smooth~functions~or~sections} \\
{\rm flat~bundles}     & \leftrightarrow & {\rm holomorphic~bundles} \\
{\rm locally~constant~functions~or~sections}
                       & \leftrightarrow &
{\rm local~holomorphic~functions~or~sections} \\
{\rm cohomology~of~locally~constant~sheaves}
                       & \leftrightarrow &
\mbox{cohomology of sheaves of ${\cal O}_X$-modules} \\
\end{array}
$$
(Here ${\cal O}_X$ denotes  the  sheaf of holomorphic functions on a complex
(algebraic) manifold $X$.) Very informally, this table could be summarized
in one line with

\nopagebreak
\centerline{``Topology" versus ``Complex Algebraic Geometry".}

Our interest in this line of thinking is related to the ideas of Arnold
on complex analytic analogues of differential
geometric concepts (cf., \cite{Arn}).
Some features of the above correspondence can also
be found in the papers \cite{FK,DT,KR}. In particular,
the approach of Donaldson and Thomas \cite{DT} of
transferring differential
geometric constructions into the context of complex analytic
(or algebraic) geometry could
lead one to a complexification of geometry in a sense similar
to  the complexification of topology pursued here.

There are also several motivations from mathematical physics:
in particular, from considering any topological field theory
of type B \cite{ASL,arhar} and of BV type \cite{Sch}
or, e.g., a complex analogue of the
Chern--Simons gauge theory suggested in Ref.~\cite{W}. The latter
context leads us immediately to a search for a proper holomorphic
analogue of the linking number (cf., also \cite{Ger, FT}).

\begin{ssect}{Holomorphic orientation.} \label{InPair}
Let $X$ be a compact complex manifold and $u$  be a smooth
$(0,k)$-form on it, $0\leqslant k\leqslant n=\dim\,X$.
We would like to treat such $(0,k)$-forms in the same manner as
ordinary $k$-forms on a smooth manifold, but in the framework
of complex geometry.  In particular, we have to be able to
integrate them  over $k$-dimensional {\it complex} submanifolds in $X$.
Recall that in the theory of differential forms, a form can be
integrated over a real submanifold provided
that the submanifold is
endowed with an orientation. Thus, we need to find a holomorphic
analogue of the orientation.

Obviously, if a $k$-dimensional submanifold $W\subset X$ is equipped
with a holomorphic
$k$-form  $\omega$, one can consider the following integral
$$
\int_W \omega\wedge u \;  
$$
of the product of the $(k,0)$- and $(0,k)$-forms.
Therefore we are going to regard a top degree holomorphic form on a complex
manifold as an analogue of orientation. More generally, if the form $\omega$
is allowed to have first order poles on a smooth hypersurface in $W$, the above
integral is still well-defined.

\end{ssect}

\begin{ssect}{The Cauchy--Stokes formula.}  \label{CSt}
The new feature brought by the presence of poles of
$\omega$ shows up in the following relation.
Consider the integral (\ref{InPair}) with a meromorphic $k$-form $\omega$
having first order poles on a smooth hypersurface $V\subset W$.
Let the smooth $(0,k)$-form $u$ on $X$ be
$\bar\partial\,$-exact, that is $u=\bar\partial\,v$ for some
$(0,k-1)$-form $v$ on $X$. Then
$$
\int_W \omega\wedge\bar\partial\,v
= 2\pi i\int_V\res\,\omega\wedge v \;.\label{CSf}
$$
We shall exploit this straightforward generalization of the Cauchy
formula as a comp\-lexi\-fied analogue of the Stokes theorem.
Here $\res\,\omega$ denotes a  $(k-1)$-form on $V$
which is the Poincar\'e residue of $\omega$ (see Sect.~\ref{PR}).
\end{ssect}

\begin{ssect}{Boundary operator.}
The formula (\ref{CSt}) prompts  us to consider
the pair $(W,\omega)$ consisting of a $k$-dimensional submanifold
$W$ equipped with a meromorphic form $\omega$
(with first order poles on $V$) as an
analogue of a compact oriented submanifold with boundary. In the
present paper we construct a homology theory in which the
pairs $(W,\omega)$ will play the role of chains, while the boundary
operator will take the form $\partial\,(W,\omega)=\tpi(V,\res\,\omega)$.
Note, that in the situation under consideration, when the polar set
$V$ of the form $\omega$ is a smooth $(k-1)$-dimensional submanifold in a
smooth $k$-dimensional $W$, the induced ``orientaion" on $V$ is given
by a regular $(k-1)$-form $\res\,\omega$. This means that
$\partial\,(V,\res\,\omega)=0$, or the boundary of a boundary is zero.
The latter will be the source of the identity $\partial^2=0$ in the homology
theory discussed below. We shall call it  the {\em polar
homology.}
\end{ssect}

\begin{ssect}{Pairing to smooth forms.} \label{InHom}
It is clear that the (would-be) polar homology groups of a complex
manifold $X$ should have a pairing to Dolbeault cohomology groups
$H^{0,k}_{\bar\partial}(X)$. Indeed, for a polar $k$-chain $(W,\omega)$
and any $(0,k)$-form $u$ such a pairing is given by the integral
$$
\langle(W,\omega)\,,u\rangle = \int_W \omega\wedge u \;.
$$
In other words, the polar chain $(W,\omega)$ defines a current on $X$ of
degree \linebreak $(n,n-k)$, where $n=\dim X$.
This pairing descends to
(co)homology classes by virtue of the Cauchy--Stokes
formula (\ref{CSt}), see Sect.~\ref{Hom}.

\end{ssect}

\begin{ex} \label{InCurv}
Now we are already able to find out what
could be the polar homology groups $HP_k$
of a complex projective curve $X$.  In this (and in any)
case, all the 0-chains are cycles. Let $(P,a)$ and $(Q,b)$ be
two 0-cycles, where
$P,Q$ are points on $X$ and $a,b\in\C\,$. They are polar
homologically equivalent iff $a=b$.  Indeed, $a=b$ is necessary and
sufficient for the existence of a meromorphic 1-form $\alpha$ on $X$,
such that $\dvsr_\infty\alpha = P+Q$ and
$\res_P\,\alpha=\tpi\,a,\, \res_Q\,\alpha=-\tpi\,b$.
(The sum of all residues of a meromorphic differential on a
projective curve is zero by the Cauchy theorem.)
Then, we can write in
terms of polar chain complex (to be defined in detail
in Sect.~\ref{PH}) that
$(P,a)-(Q,a)=\partial\,(X,\alpha)$. Thus, $HP_0(X)=\C\,$.

As to polar 1-cycles, these correspond to all possible holomorphic
1-forms on $X$. On the other hand, there are no 1-boundaries,
since there are no polar 2-chains in $X$. Hence $HP_1(X)\cong\C\,^g$, where
$g$ is the genus of the curve $X$.
(In particular, the polar Euler characteristics of $X$ equals $1-g$
and coincides with its holomorphic Euler characteristics.)

Similar considerations show that  for any $n$-dimensional $X$ we have
$HP_n(X)\, =\, H^0(X,\Omega^n_X)$ and, if $X$ is connected, also
$HP_0(X)=\C\,$.
\end{ex}

\begin{ssect}{Polar intersections.}
One can define a complex (polar) analogue  of the intersection number
in topology. For instance, let
$(X, \mu)$ be a complex manifold equipped with
a meromorphic volume form $\mu$ without zeros (its ``polar orientation").
Consider two polar cycles
$(A,\alpha)$ and $(B,\beta)$ of complimentary dimensions that
intersect transversely in $X$ (here $\alpha$ and $\beta$ are
volume forms, or ``polar orientations," on the corresponding submanifolds).
Then the polar  intersection number is defined by the formula
$$
        \langle(A,\alpha) \Cdot (B,\beta)\rangle =
    \sum_{P\in A\cap B}\frac{\alpha(P)\wedge\beta(P)}{\mu(P)} \;\,.
$$
(For explanations, see (\ref{Xtrans}).)
At every intersection point $P$, the ratio in the right-hand-side
is the ``comparison" of the orientations of the polar cycles at that
point with  the orientation of the ambient manifold.
This is a straightforward analogue of the use of mutual orientation
of cycles in the definition of the topological intersection number.
Note, that in the polar  case the intersection number does not have
to be an integer.
(Rather, it is a holomorphic function of the ``parameters"
$(A,\alpha), (B,\beta)$ and $(X,\mu)$.)

Similarly, there is a polar analogue of the
intersection product of cycles when they intersect over
a manifold of positive dimension (see Sect.~\ref{X}).

\end{ssect}

\begin{ssect}{Polar links.}
By developing this analogy further we come to a polar analogue of the
linking number. For instance, in complex dimension three we start with
two smooth polar 1-cycles $(C_1,\alpha_1)$ and $(C_2,\alpha_2)$,
i.e.\ $C_1$ and $C_2$ are smooth complex curves equipped with
holomorphic 1-forms
in a three dimensional $X$.
Let us take the 1-cycles which are polar boundaries. This means, in
particular, that there exists such a 2-chain $(S_2,\beta_2)$  that
$(C_2,\alpha_2)=\partial\,(S_2,\beta_2)$.
Suppose, the curves $C_1$ and $C_2$ have no common points and
$S_2$ is a smooth surface which
intersects transversely with the curve $C_1$. Then, in analogy with the
topological linking number of two curves in a three-fold, we define
the polar linking number of the 1-cycles above
as the polar intersection number of the 2-chain  $(S_2,\beta_2)$
with the 1-cycle $(C_1,\alpha_1)$:
$$
\lk\left( (C_1,\alpha_1), (C_2,\alpha_2)\right):=
\sum\limits_{P\in C_1\cap \,S_2} \frac{\alpha_1(P)\wedge\beta_2(P)}{\mu(P)}\,.
$$
One can show that the expression above does not depend on the choice of
$(S_2,\beta_2)$,  and has certain invariance
properties mimicking those of the topological linking number
in ``polar'' language.
We are going to discuss the properties of $\lk$ in more detail
in a future publication.

\end{ssect}

\begin{rem}
Most of the above  discussion extends to polar chains
$(A,\alpha)$  where the meromorphic $p$-form $\alpha$ is not
necessarily of top degree, that is $0\leqslant p\leqslant k$,
where $k=\dim_{\C} A$.
To define the boundary operator we have to restrict ourselves
to the meromorphic forms with logarithmic singularities.
The corresponding
polar homology groups enumerated by two indices $k$ and $p$
($0\leqslant p\leqslant k$) will
be discussed elsewhere (see, though, some remarks in Sect.~\ref{Hom}.b below).
\end{rem}


\mysection{Preliminaries} \label{Res}
\nopagebreak

{\sl a) Polar divisors and  residues}
\nopagebreak
\medskip
\nopagebreak

\noindent
The Poincar\'e residue is a higher-dimensional generalization
of the classical Cauchy residue,
where the residue at a point in a domain of one complex variable
is generalized to the residue at a hypersurface.
\smallskip

\begin{ssect}{} \label{PR}
Let $M$ be an $n$-dimensional complex manifold
and $\omega$ be a meromorphic $n$-form on $M$ which is allowed
to have  first order poles on a smooth hypersurface $V$.
Then, the form $\omega$ can be locally expressed as
$$
\omega = \frac{\rho\wedge dz}{z} + \varepsilon\;,
$$
where $z=0$ is a local equation of $V$ and $\rho$ (respectively,
$\varepsilon$) is a
holomorphic $(n-1)$-form (resp., $n$-form).
Then the restriction $\rho|_V$ is an unambiguously defined
holomorphic $(n-1)$-form on $V$.
\end{ssect}

\begin{defn} \label{PRdef}
The Poincar\'e residue of the $n$-form $\omega$ in {\rm(\ref{PR})} is the
following $(n-1)$-form on~~$V$
$$
\res\,\omega := \rho|_V \; .
$$
\end{defn}

\smallskip
\begin{ssect}{} \label{rep-res}
It is straightforward to extend this to the case of normal crossing
divisors. Suppose that the meromorphic
$n$-form $\omega$ has the first order poles on a normal crossing
divisor $V=\cup_i V_i$ in $M$. [Normal crossing divisor means
that $V$ has only smooth components $V_i$ (each entering with
multiplicity one) that intersect generically.]
Analogously to the Definition~\ref{PRdef} one can define a residue at each
component $V_i$. The resulting $(n-1)$-forms $\res_{V_j}\omega$ are then
meromorphic and  have first order poles at
the pairwise intersections $V_{ij}=V_i\cap V_j$
One can now consider
the repeated Poincar\'e residue at $V_{ij}$.
Representing $\omega$ as
$\omega=\varrho\wedge\frac{dz_i}{z_i}\wedge\frac{dz_j}{z_j}$,
where $z_i=0$ and $z_j=0$ are local equations of the components
$V_i$ and $V_j$ respectively one finds that
$$
\res_{i,j}\omega:=\res_{V_{ij}}\left(\res_{V_j}\omega\right)
=\res_{z_i=0}\left(\res_{z_j=0}\,\varrho\wedge\frac{dz_i}{z_i}
\wedge\frac{dz_j}{z_j}\right)=\varrho|_{V_{ij}}.
$$
Note that
$$
\res_{i,j}\omega=-\res_{j,i}\omega\,.
$$

\smallskip

{\sc Notation.} Let us denote by $\res\,\omega$ the collection of
$(k-1)$-forms $\res_{V_j}\omega$, the residues of $\omega$ at the components
of the normal crossing divisor $\dvsr_\infty\omega = \cup_i V_i$.

\end{ssect}

\medskip

{\sl b) The push-forward map} (see \cite{Gr})

\medskip

\noindent
For a finite covering $f\!:X\to Y$ and a function $\varphi$ on $X$ one
can define its push-forward, or the trace, $f_*\,\varphi$ as a function
on $Y$ whose value at a point is calculated by summing over the
preimages taken with multiplicities. The operation $f_*$ can be
generalized to $p$-forms and to the maps $f$ which are only
generically finite.
\smallskip

\begin{ssect}{} \label{pu-fo}
Suppose that $f\!:X\to Y$ is
a proper, surjective holomorphic mapping
where both $X$ and $Y$ are smooth complex manifolds of the same
dimension $n$. The {\it push-forward map} is a mapping
$$
f_*:\Gamma(X,\Omega^p_X)\to \Gamma(Y,\Omega^p_Y)\;.
$$
The push-forward map is also defined  for meromorphic forms,
$ f_*:\Gamma(X,{\cal M}^p)\to \Gamma(Y,{\cal M}^p) $.

Its construction is as follows.
First note that $f$ is generically finite, i.e., there is an analytic
hypersurface $D\subset Y$ such that $f$ is finite unramified covering
away from this hypersurface $D$. Hence, for sufficiently small
open neighborhood $U$ of any point in $Y^*:=Y\smallsetminus D$, the
inverse image
$f^{-1}(U)=U_1\sqcup\ldots\sqcup U_d$
is a disjoint union of $d$ open sets $U_j$, such that
$f|_{U_j}$ is an isomorphism with the inverse $s_j: U\to U_j$.
Given a form $\omega$ on $X$, one defines its push-forward
$$
f_*\,\omega:=s^*_1\,\omega+...+s^*_d\,\omega
$$
in $U$, and therefore, in $Y^*$.
One can check that the form $f_*\,\omega$
extends across the smooth points of $D$ and, hence, to the whole of
the manifold $Y$,
since the remaining part of $D$ has codimension greater than one.
The resulting form $f_*\,\omega$
is holomorphic (resp. meromorphic) on $Y$
provided the form $\omega$ was holomorphic (resp. meromorphic) on $X$.
\end{ssect}

The operations of push-forward and residue are related in the following way.

\begin{prop} \label{p-f-res}
Let $\pi\!:X\to Y$ be a proper, surjective holomorphic map of complex
manifolds of the same dimension $n$. Let $\omega$ be a meromorphic form with
only first order poles on a smooth hypersurface $V$ in $X$. Suppose
$V_o:=\pi(V)$ is a smooth hypersurface in $Y$. Then $\pi_*\omega$ has first
order poles on $V_o$ and
$$
\res\,\pi_*\omega = \bar\pi_*\,\res\,\omega \;,
$$
where $\bar\pi\!:V\to V_o$ is the restriction of $\pi$.
\end{prop}



\mysection{Polar Homology of Projective Varieties} \label{PH}

Here we define a homological complex based on the notion of the polar
boundary.  The construction is analogous to the definition of
homology of a topological space with replacement of continuous maps by
complex analytic ones. The notion of the boundary (of a simplex or a
cell) is replaced by the Poincar\'e residue of a meromorphic
differential form. There are however
important distinctions. First, we shall only have an analogue of
the non-torsion part of
homology. Second, unlike the topological
homology, where in each dimension $k$ one uses all continuous maps of
one standard object (the standard $k$-simplex or the standard $k$-cell)
to a given topological space, in polar homology we deal with complex
analytic maps of a large class of $k$-dimensional varieties to a given
one.

\begin{ssect}{Polar chains.}\label{smooth}
In this section we deal with complex projective varieties,
i.e., subvarieties of a complex projective space.
(In this setting
the complex analytic considerations are equivalent to
algebraic ones.)
By a smooth projective variety we always understand a smooth and connected
one. For a smooth variety $M$, we denote by
$\Omega^p_M$ the sheaf of regular $p$-forms on $M$.
The sheaf $\Omega^{\dim M}_M\,$ of forms of
the top degree on $M$ will sometimes be denoted by $K_M$.

         The space of polar $k$-chains for a complex
projective variety $X, \dim X=n,$
will be defined as a $\C\,$-vector space with certain generators and
relations.
\end{ssect}

\begin{defn} \label{PHdef}
The space of {\sl polar $k$-chains~} ${\cal C}_k(X)$ is
a vector space over $\C$ defined as the quotient
${\cal C}_k(X)=\hat{\cal C}_k(X)/{\cal R}_k$, where the vector space
$\hat{\cal C}_k(X)$ is freely generated by the
triples $(A,f,\alpha)$ described in (i),(ii),(iii) and
${\cal R}_k$ is defined as relations
{\rm (R1),(R2),(R3)} imposed on the triples.

(i) $A$ is a smooth complex projective variety, $\dim A=k$;

(ii) $f\!: A\to X$ is a holomorphic map of projective varieties;

(iii) $\alpha$ is a rational $k$-form on $A$
with first order poles on $V\subset A$, where $V$ is a normal crossing
divisor in $A$, i.e., $\alpha\in\Gamma(A,\Omega^k_A(V))$.

\noindent
The relations are:

{\rm (R1)} $\lambda (A, f, \alpha)=(A, f, \lambda\alpha)$

{\rm (R2)} $\sum_i(A_i,f_i,\alpha_i)=0$ provided that
 $\sum_if_{i*}\alpha_i\equiv 0$, where
$\dim f_i(A_i)=k$ for all $i$ and
the push-forwards $f_{i*}\alpha_i$ are considered on the smooth part
of $\cup_i f_i(A_i)$;

{\rm (R3)} $(A,f,\alpha)=0$ if $\dim f(A)<k$.

\end{defn}

{\sl Remarks to the definition}
\medskip

\begin{ssect}{}
By definition, ${\cal C}_k(X)=0$ for $k<0$ and $k>\dim X$.

\end{ssect}

\begin{ssect}{}
In what follows we sometimes will make no difference between a triple
$(A,f,\alpha)$ and the equivalence class defined by it in ${\cal C}_k(X)$.
An arbitrary polar chain can thus be written as a sum of triples of the form
$\sum_i(A_i,f_i,\alpha_i)$. A chain equivalent to a single triple will be
called {\em prime}.
For a chain $a=\sum_i(A_i,f_i,\alpha_i)$, let us call
the subvariety $\cup_if_i(A_i)$ in $X$ the {\em support} of $a$. If the
support of a chain is a smooth subvariety in $X$, such a chain will be
called {\em smooth}. One can show that smooth chains are prime,
since we suppose that ``smooth" implies ``connected" (see \ref{smooth}).

\end{ssect}

\begin{ssect}{}
The relation (R2) allows us, in particular, to refer to prime polar
chains as pairs replacing a triple $(A,f,\alpha)$ by a pair
$(\hat A,\hat\alpha)$, where $\hat A=f(A)\subset X$, $\hat\alpha$
is defined only on
the smooth part of $\hat A$ and $\hat \alpha=f_*\alpha$ there.
Due to the relation (R2), such a pair $(\hat A,\hat\alpha)$ carries
precisely
the same information as $(A,f,\alpha)$.\footnote{Note, that the
consideration of triples $(A,f, \alpha)$
instead of pairs $(\hat A,\hat\alpha)$ which we used in Sect. 1
is similar to
the definition of chains in the singular homology theory:
in the latter case, although one considers
the mappings of abstract simplices into the manifold,
but morally
it is only ``images of simplices'' that matter.  }
(The only point to worry about
is that such pairs cannot be arbitrary. In fact, by the Hironaka
theorem on resolution of singularities, any subvariety $\hat A\subset X$
can be the image of some regular $A$, but the form $\hat\alpha$ on the
smooth part of $\hat A$ cannot be arbitrary.)

\end{ssect}

\begin{ssect}{}
The relation (R2) also represents additivity with respect to $\alpha$,
that is
$$
(A,f,\alpha_1)+(A,f,\alpha_2)=(A,f,\alpha_1+\alpha_2).
$$
Formally speaking, the right hand side makes sense only if
$\alpha_1+\alpha_2$ is an admissible form on $A$, that is if its polar
divisor $\dvsr_\infty(\alpha_1+\alpha_2)$ has normal crossings.
However, one can always replace $A$ with a variety $\tilde A$,
obtained from $A$ by a blow-up,
$\pi\!:\tilde A\to A$, in such a way
that $\pi^*(\alpha_1+\alpha_2)$ is admissible on $\tilde A$, i.e.,
$\dvsr_\infty(\alpha_1+\alpha_2)$ is already a normal crossing
divisor. (This is again the Hironaka theorem.) The (R2) says that
$(A,f,\alpha_1)+(A,f,\alpha_2)=(\tilde A,f\circ\pi,\pi^*(\alpha_1+\alpha_2))$.

\end{ssect}

\begin{defn} \label{d-def}
The {\sl boundary operator~} $\partial: {\cal C}_k (X)\to{\cal C}_{k-1}(X)$
is defined by
$$
\partial(A,f,\alpha)=\tpi\sum_i(V_i, f_i, \res_{V_i}\,\alpha)
$$
(and by linearity),
where $V_i$ are the
components of the polar divisor of $\alpha$,
$\dvsr_\infty\alpha=\cup_iV_i$, and the maps $f_i=f|_{V_i}$
are restrictions of the map $f$ to each component of the divisor.

\end{defn}

\begin{th}
The boundary operator $\partial$ is well defined, i.e.\ it is
compatible with the relations {\rm (R1),(R2),(R3)}.

\end{th}

\begin{proof}{} We have to show that $\partial$ respects the relations
(R1),(R2),(R3), in other words, $\partial$ maps equivalent
sums of triples to equivalent ones. It is trivial with (R1). To check (R2),
let us recall Proposition~\ref{p-f-res}.
Consider a sum of triples $\sum_i (A_i,f_i,\alpha_i)$
belonging to (R2), that is
$\dim A_i=\dim f_i(A_i)=k, \forall i$, and
$\sum_i f_{i*}\alpha_i=0$ on the smooth part of $\cup_i f_i(A_i)$.
Since the irreducible components of $\cup_i f_i(A_i)$ can be treated
separately it is natural to consider only the case when all the
triples have the same support, $f_i(A_i)=\hat A \subset X, \forall i$. Let
$V_i\subset A_i$ be the divisor of poles of $\alpha_i$ and let
$\hat V:=\cup_i f_i(V_i)\subset \hat A$. We want to prove that
$$
\sum_i \bar f_{i*}\res\,\alpha_i=0
$$
on the smooth part of $\hat V$, where
$\bar f_{i*}\,: V_i\to\hat V$  for each $i$ is the restriction
of the map $f_i$.

Suppose first that there exists a smooth point of $\hat V$
which is smooth also
in $\hat A$. Then the Proposition~\ref{p-f-res} applied in a neighborhood of
that point gives us the desired vanishing
$\sum_i \bar f_{i*}\res\,\alpha_i=0$, as a consequence of the
equality $\sum_i f_{i*}\alpha_i=0$. This is however not enough for our proof
since some components of $\hat V$ may lie entirely in the set of singular
points of $\hat A$. To overcome this problem we apply the Hironaka theorem
replacing $\hat A$ with a smooth variety $\tilde A$, a blow-up of $\hat A$,
and correspondingly blowing up all $A_i$, so that the following diagram is
commutative:
$$
\begin{array}{lrl}
\phantom{\mbox{\footnotesize $f_i$}}A_i~
& \stackrel{\pi_i}{\longleftarrow}
& \phantom{\mbox{\footnotesize $\tilde f_i$}}\tilde A_i  \\
\mbox{\footnotesize $f_i$}\downarrow~
&
& \mbox{\footnotesize $\tilde f_i$}\downarrow   \\
\phantom{\mbox{\footnotesize $f_i$}}\hat A~
& \stackrel{\pi}{\longleftarrow}
& \phantom{\mbox{\footnotesize $\tilde f_i$}}\tilde A
\end{array}
$$
Then we apply Proposition \ref{p-f-res} on the blown up side.

We must recall now that the divisor
$V_i=\dvsr_\infty\alpha_i$ could have components
that were mapped by $f_i$ to subvarieties of dimension less than
$k-1$; hence, we conclude that we have just proved the following statement
(symbolically): if $a\in\,$(R2) then $\partial a\in\,$(R2)+(R3).

Now, it remains to prove the compatibility of $\partial$ with (R3). Let
$a=(A,f,\alpha)$ be a degenerate triple described in (R3), i.e.,
$\dim f(A)<k=\dim A$. We shall show that $\partial a\in\,$(R2)+(R3)
in this case. The polar divisor $V=\dvsr_\infty\alpha,~\dim V=k-1$, is,
by assumptions of Definition~\ref{PHdef}, a normal crossing divisor in $A$.
Let us split the components of $V$ into two parts: non-degenerate and
degenerate ones. That is $V=N\cup D$ where $\dim f(N)=k-1$ and
$\dim f(D)<k-1$. According to this splitting, $\partial a$ is represented
as a sum of two terms corresponding to $\res_N\alpha$ and $\res_D\alpha$.
The second term belongs to (R3) and we have to show only that
the first one belongs to (R2), i.e., that
$\bar f_*\res_N\alpha=0$, where $\bar f=f|_N$. Recall that we
suppose that $\dim f(A)<k$. If it happens that $\dim f(A)<k-1$,
we have $N=\varnothing$ and there is nothing to prove. Therefore we may
assume that $\dim f(A)=k-1$ and, by irreducibility of $A$, $f(A)=f(N)$.
Then, for a generic smooth point $Q\in f(A)$, its preimage in $A$,
$C:=f^{-1}(Q)\subset A$, is a smooth projective curve. This curve intersects
with $N$ over the set $\bar f^{-1}(Q)$ and we may suppose that the latter
consists of a finite number of points $P_i$ which are smooth in $N$ and that
the intersections are transverse there.

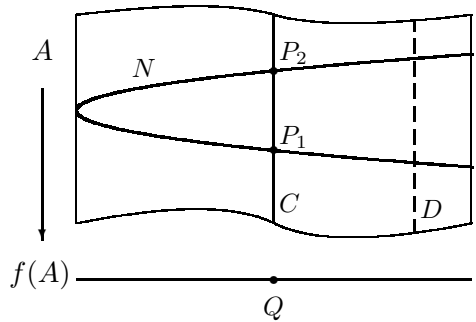
\myfigure

Let $\beta(P_i)$ denote the values
of $\res_N\alpha$ at the points $P_i\in N\cap C$ and pick up a non-vanishing
$(k-1)$-form $\beta_o$ at $Q$ (recall that $Q=f(P_i),~\forall i\,$). Let us
show that
$$
\sum_i \frac{\beta(P_i)}{f^*\beta_o(P_i)} = 0
$$
(this would mean that $\bar f_*\res_N\alpha=0$ on the smooth part of $\bar
f(N)=f(A)$ --- the required result).
To prove this, let us notice that there
exists a meromorphic 1-differential $\omega$ on $C$ such that
$$
\omega(P)\otimes f^*\beta_o(P) = \alpha(P),~~P\in C \;.
$$
($\omega$ is obtained by dividing $\alpha$ by the non-vanishing form
$f^*\beta_o\,$.)
This equality is understood in the sense of the natural isomorphism
$$
K_C\otimes\left.\left(f^*K_{f(A)^*}\right)\right|_C =
\left.K_A\right|_C  \;,
$$
where $f(A)^*$ is the smooth part of $f(A)$.
It is easy to see now that for $\beta(P_i)=\res_N\alpha(P_i)$ we have
$$
\sum_i \frac{\beta(P_i)}{f^*\beta_o(P_i)} =
\sum_i \res_{P_i}\omega = 0  \;.
$$
The latter equality follows from the observation that $P_i$ are the
only points on $C$ where $\omega$ has poles. Indeed, the poles of $\omega$
are located on
$\dvsr_\infty\alpha\cap C=(N\cap C)\cup(D\cap C)$.
One part of this gives us the points $P_i,~\{P_i\}=N\cap C$, while the
rest, $D\cap C$,
corresponding to the ``degenerate" part $D$ of $\dvsr_\infty\alpha$
can be assumed to be empty,
$D\cap C=\varnothing$. Indeed, we could have assumed from the very
beginning that
$C=f^{-1}(Q)$ does not meet $D$ because $\dim f(D)<k-1$ and we might suppose that
$Q\notin f(D)$.
\end{proof}

\begin{th}
~~~$\partial^2=0\;$.
\end{th}

\begin{proof}{} We need to prove this for triples
$(A,f,\alpha)\in {\cal C}_k(X)$, i.e.,
for forms $\alpha$ with normal crossing divisors of poles.
The repeated residue at pairwise intersections differs by a sign according to
the order in which the residues are taken, see \ref{rep-res}.
Thus the contributions to the repeated residue from different
components cancel out (or, the residue of a residue is zero).\footnote
{One can note that
an example of the polar divisor $\{xy=0\}$
for the form $dx\wedge dy/xy$ in ${\Bbb C}^2$
should be viewed as a complexification of a polygon vertex in
${\Bbb R}^2$. Indeed,  the cancellation of the repeated residues
on different components of the divisor is mimicking the calculation of the
boundary of a boundary of a polygon: every polygon vertex
appears twice with different signs as a boundary point of two sides.}
\end{proof}

\begin{defn}
For a smooth complex projective variety $X, \dim X=n$, the chain
complex
$$
0 \to {\cal C}_n(X) \stackrel{\!\!\partial}{\longrightarrow}
{\cal C}_{n-1}(X) \stackrel{\!\!\partial}{\longrightarrow} \dots
\stackrel{\!\!\partial}{\longrightarrow} {\cal C}_0(X)\to 0
$$
is called the {\sl polar chain complex} of $X$. Its homology groups,
$HP_k(X), k=0,\ldots,n$, are called the {\sl polar homology groups} of $X$.

\end{defn}

\begin{ex} For a projective curve of genus $g$ the polar homology
groups are as follows: $HP_0=\C\,, \, HP_1=\C^{\,g}$,
and $HP_k=0$ for $k\geq 2$.
One can readily see that the approach with triples coincides
with the consideration of Introduction.
\end{ex}

\begin{rem}
The functorial properties of polar homology are straightforward.
A regular morphism of projective varieties
$h\!:X\to Y$ defines a homomorphism
$h_*\!:PH_k(X)\to PH_k(Y)$.

\end{rem}

\begin{rem} \label{non-top}
The definitions of polar chains can be generalized to the case of
$p$-forms on $k$-manifolds, i.e., for the forms of
not necessarily top degree, $p\leq k$.
Instead of meromorphic $k$-forms with poles of the first order we
have to restrict ourselves by $p$-forms with logarithmic
singularities. The definition of the boundary operator $\partial$,
the  property $\partial^2=0$, and the definition of the polar homology groups
can be carried over to this, more general, situation. The polar
homology groups
are then enumerated by two indices: $HP_{k,p}(M)$. The definition above
corresponds to the $p=k$ case.
We will discuss the more general polar homology groups elsewhere [iKR].
\end{rem}

\begin{ssect}{Relative polar homology.} \label{relPHdef}
Let $Z$ be a projective subvariety in a projective $X$.
Analogously to the topological relative homology we can define the
polar relative homology of the pair $Z\subset X$.

\end{ssect}

\begin{defn} The relative polar homology $HP_k(X,Z)$ is the
homology of the following quotient complex of chains:
$$
{\cal C}_k(X,Z)\:={\cal C}_k(X)/{\cal C}_k(Z).
$$
{\rm
Here we use the natural embedding of the chain groups
${\cal C}_k(Z)\hookrightarrow{\cal C}_k(X)$.
This leads to the long exact sequence in polar homology:
$$
\ldots\to HP_k(X)\to HP_k(X,Z) \stackrel{\!\!\partial}{\longrightarrow}
HP_{k-1}(Z)\to HP_{k-1}(X)\to\ldots
$$
}
\end{defn}

\begin{ssect}{Systems of coefficients.}  \label{coeff}
One can introduce the notion of a homological system of coefficients
appropriate for the polar complex. The most geometrical example would
be, perhaps, to supply projective varieties $A,~ f\!:A\to X,$
with coherent sheaves
${\cal F}_{A,f}$ obeying certain relations between
${\cal F}_{A_1,f_1}$ and ${\cal F}_{A_2,f_2}$ when $f_1(A_1)=f_2(A_2)$ and related
by some homomorphisms playing the role of the residue. We do not
study this in the present paper, but let us mention that the homology
groups appearing in  Sects.~\ref{non-top} above and
\ref{other} below can be viewed as an example. On
the other hand, the simplest case of a polar homological system of
coefficients corresponds to
${\cal F}_{A,f}=f^*{\cal F}\otimes K_A$,
where ${\cal F}$ is a sheaf on $X$ and $\alpha$ in the triple
$(A,f,\alpha)$ is understood as a global section of
$f^*{\cal F}\otimes K_A(V)$.
Let us denote the corresponding
homology as $HP_k(X,{\cal F})$. This case is mentioned in
Sects.~\ref{Conj}, \ref{curve}.

\end{ssect}


\mysection{Polar Chains and Differential Forms}
\label{Hom}

{\sl a) Dolbeault cohomology as polar de Rham cohomology}

\medskip

\noindent
As we discussed in the Introduction,
the Dolbeault complex of $(0,k)$-forms should be related
to the polar homology in the same  way as
the de Rham complex of smooth forms is
related to the topological homology (e.g., singular homology).
Now, after the definitions of Sect.~\ref{PH} are given,
we are able to make this point more explicit.

\smallskip

\begin{ssect}{}\label{pair}
In a smooth projective variety $X$, consider a polar $k$-chain,
for instance, a prime one,
i.e.\  (an equivalence class of) a triple $a=(A, f, \alpha)$.
Such a triple can be regarded as a linear functional on the space
of smooth $(0,k)$-forms on $X$.
Let $u$ be a smooth $(0,k)$-form on $X$, then the pairing is given by the
following integral:
$$
\langle a\,,u\rangle\: := \int_A\alpha\wedge f^*u \;.
$$
The integral is well defined since $\alpha$
has only first order poles on a normal crossing divisor.
It is now straightforward to show that the pairing $\langle\;,\,\rangle$
descends to the space of equivalence classes of triples ${\cal C}_k(X)$,
i.e., that it is compatible with the relations (R1), (R2), (R3) of
Definition~\ref{PHdef}. Indeed, (R1) is obvious, compatibility with (R3)
follows from noticing that $f^*u=0$ if $\dim f(A)<k$, and the compatibility
with (R2) follows from the relation
$\int_A\alpha\wedge f^*u=\int_{f(A)}f_*\alpha\wedge u$ if $\dim f(A)=k$,
where the last integral is taken over the smooth part of $f(A)$.

\end{ssect}

\begin{rem} \label{currents}
Let us notice that the last considerations can be
used\footnote{We thank A.~Levin for emphasizing this point.}
as an alternative definition of the polar chain complex on a smooth
projective variety $X$ (or any smooth compact complex manifold).
The pairing above can be thought of as a map
$\hat\varphi\,:\hat{\cal C}_k(X)\to{\cal D}^{n,n-k}(X)$,
where $\hat{\cal C}_k(X)$ is the vector space freely generated by the
triples $(A,f,\alpha)$ (see Definition \ref{PHdef}) and
${\cal D}^{n,n-k}(X)$ is the space of currents of degree $(n,n-k)$ on $X$
which is defined as a space of certain linear functionals on the smooth
$(0,k)$-forms (see \cite{GH}). Then the relations (R1),(R2),(R3) in the
Definition~\ref{PHdef} correspond to the kernel of the map $\hat\varphi$. In
other words, the space of polar chains ${\cal C}_k(X)$ can be defined as a
subspace of currents --- the image of $\hat\varphi$. We have thus an
embedding
$$
\varphi:\, {\cal C}_k(X) \hookrightarrow {\cal D}^{n,n-k}(X) \,.
$$
Moreover, the Cauchy--Stokes formula (\ref{CSt}) shows that
${\cal C}_k(X)$ is a subcomplex of the $\bar\partial$-complex of currents
${\cal D}^{n,n-k}(X)$, i.e.\ for $a\in{\cal C}_k(X)$ we have
$\varphi(\partial a)=\bar\partial\varphi(a)$. (This is in fact shown also in
the proof of \ref{hom} below.)

\end{rem}

\begin{prop} \label{hom}
The pairing {\rm(\ref{pair})} defines the following
homomorphism in (co)homology:
$$
\rho\!: HP_{k}(X)\to H^{n,n-k}_{\bar\partial}(X) \;,
$$
where $n=\dim X$.
\end{prop}

\begin{proof}{} By the Serre duality, $\rho$ is
the map  $HP_{k}(X)\to ( H^{0,k}_{\bar\partial}(X))^*$ ~and it
is sufficient to verify that the
pairing (\ref{pair}) vanishes if
$\partial\,a=0$ and $u=\bar\partial\, v$,
or if $\bar\partial\, u=0$ and $a=\partial\,b$.
This follows
immediately from the Cauchy--Stokes formula (\ref{CSt}):
$$
\int\limits_A\alpha\wedge f^*(\bar\partial u)
\:=\: \tpi\!\! \int\limits_{\dvsr_\infty\alpha}\!(\res~\alpha)\wedge f^*(u) \;,
$$
that is $\langle a\,,\bar\partial u\rangle =
\langle\partial a\,,u\rangle$.
\end{proof}

A number of examples suggests that, for projective manifolds,
the homomorphism (\ref{hom}) should be in fact an isomorphism.

\begin{conj} \label{Conj}
{\bf (Polar de Rham theorem)}
For a smooth projective manifold $X$ the mapping
$\rho: HP_{k}(X)\to H^{n,n-k}_{\bar\partial}(X)$ is an isomorphism of
the polar homology and Dolbeault cohomology groups. Equivalently, in
terms of dual cohomology groups,
$$
HP^{k}(X)\cong  H^{0,k}_{\bar\partial}(X) \;.
$$
\end{conj}

\noindent
An analogous conjecture that
$HP_k(X,{\cal F})\cong H^{n-k}(X,K_X\otimes{\cal F})$
sounds reasonable also for polar
(co)homology with coefficients in locally free sheaves on $X$
(see Remark~\ref{coeff}).

\begin{ex} \label{curve}
If $X$ is a complex curve of genus $g$
one has $HP_{0}(X)\cong\C\cong H^{1,1}_{\bar\partial}(X)$ and
$HP_{1}(X)\cong\C^{\,g}\cong H^{1,0}_{\bar\partial}(X)$
(see Example~\ref{InCurv}).
Considering also $HP_k(X,{\cal F})$ in case of ${\cal F=T}_X$ (the
tangent sheaf of $X$) one can check that
$HP_0(X,{\cal T}_X)\cong H^{0,1}_{\bar\partial}(X)$ and
$HP_1(X,{\cal T}_X)\cong H^{0,0}_{\bar\partial}(X)$.

\end{ex}

\begin{rem}
Consider the  {\it polar Euler characteristic},
$$
\chi_{pol}(X):=\sum_{k=0}^n (-1)^k\dim HP_{k}(X) \;,
$$
of an $n$-dimensional variety $X$. Then, if the conjecture
(\ref{Conj}) is true,  for a smooth projective $X$ one obtains the equality
$\chi_{pol}(X)~=~\chi_{hol}(X)$
of the polar and holomorphic Euler characteristics,
where
$$
\chi_{hol}(X)= \sum_{k=0}^n (-1)^k\dim H^{0,k}_{\bar\partial}(X)\;.
$$

\end{rem}

\bigskip


{\sl b) Forms of any degree}

\medskip

\begin{ssect}{}\label{other}
So far we considered polar chains with complex volume forms.
More generally, one could consider  polar $(k,p)$-chains $(A,f,\alpha)$, where
$\alpha$ is a meromorphic $p$-form of not necessarily maximal degree,
$p\leqslant k$, on $A$ that can have only logarithmic singularities on a
normal crossing divisor.\footnote{An important property of such forms
on projective varieties is that they are closed, see \cite{De}.}
The requirement of log-singularities is needed to have a convenient
definition of the residue and, hence, the boundary operator.

The Cauchy--Stokes formula (\ref{CSt})
extends to this case as well. As a
consequence, the natural pairing between polar $(k,p)$-chains and
smooth $(k-p,k)$-forms on $X$ gives us as before the homomorphism
(cf.\ (\ref{hom}))
$$
\rho\!: HP_{k,p}(X)\to H_{\bar\partial}^{n-k+p,n-k}(X) \;.
$$
However, unlike the case $p=k$, the map $\rho$ is not, in general,
an isomorphism for other values of $p$, $0\leqslant p<k$. For instance,
at least in the case of $p=0$, this is easy to see for the following
reason.

\end{ssect}

\begin{ssect}{Polar chains with $\mathbf p=0$.}\label{p=0}
In this case the triples $(A,f,\alpha)$ involve
0-forms (i.e., just functions) $\alpha$ on projective varieties $A$.
Then the requirement of log-singularities amounts here to saying that
these functions are holomorphic on $A$. Since $A$ is
compact, these functions must be constant. In particular,
we conclude that all
polar $(k,0)$-chains are polar cycles.

Thus, the space of polar $(k,0)$-cycles in a projective manifold $X$
is the same as the vector space generated over $\C$ by all
$k$-dimensional algebraic cycles in $X$.
(Note that the replacement of
the triples $(A,f,\alpha)$ by the
pairs $(A,\alpha)$  with $A\subset X$ is
especially convenient when $\alpha$'s are 0-forms.)
In this case one can show that the homomorphism $\rho$ maps
$HP_{k,0}(X)$ to the algebraic part of $H^{r,r}(X),$ where $ r=n-k$,
or more precisely, to
$$
H^{r,r}_{alg}(X,\C) :=
(H_{\bar\partial}^{r,r}(X)\cap H^{2r}(X,{\Q}))\otimes\C
\subset H_{\bar\partial}^{r,r}(X)  \;.
$$
This allows us to conclude that $\rho$ is not
surjective, in general. Indeed, there are examples where
$H^{r,r}_{alg}(X,\C)$ is strictly smaller than
$H_{\bar\partial}^{r,r}(X)$.  For instance, for a generic
algebraic $K3$ surface one has that $\dim H^{1,1}(X)=20$, while
$\dim H_{alg}^{1,1}(X,\C)=1$, see \cite{Tju}.
We also note that by the Hodge conjecture, the image of $\rho$
coincides with $H^{r,r}_{alg}(X,\C)$.

\end{ssect}

\begin{rem}
It would be, certainly, very interesting to describe the polar homology
groups $HP_{k,p}(X)$ for all values of $p$. In particular, it is not
clear whether the groups $HP_{k,p}(X)$ are
finite-dimensional.\footnote{It should
be mentioned, that there is a map of the complex of
polar chains to the complex considered in \cite{BO}.
The corresponding (co)homology groups are, however,
quite different, as already the simplest example of a complex curve shows.
We are grateful to S.~Bloch for illuminating discussions
on this relation.}
\end{rem}


\mysection{Intersection in Polar Homology} \label{X}
\medskip

We define here a polar analogue of the topological
intersection product.
In particular, for polar cycles of complimentary dimensions
one obtains a complex number, called the polar intersection
number.

Recall that in topology, one considers a
smooth oriented closed manifold $M$ and two oriented closed
submanifolds $A,B\subset M$ of complementary dimensions, i.e.,
$\dim A+\dim B=\dim M$. Suppose, $A$ and $B$
intersect transversely at a finite set of points.
Then to each intersection point $P$ one assigns $\pm 1$ (local
intersection index) by comparing the mutual
orientations of the tangent vector spaces $T_P A, T_P B$, and $T_P M$.

\smallskip

\begin{ssect}{Polar oriented manifolds.}
Let now $M$ be a compact complex manifold of dimension $n$, on which we
would like to define a polar intersection theory. It has to be polar oriented,
i.e., equipped with a complex volume form.
As the discussion below shows, the $n$-form $\mu$ defining its polar
orientation has to have no zeros on $M$, since we are going to consider
expressions in which $\mu$, the orientation of the ambient manifold,
enters a denominator. Therefore we
adopt the following terminology.
\end{ssect}

\begin{defn} \label{boundary-def}
i) A compact complex manifold $M$, endowed with a nowhere
vanishing holomorphic volume form $\mu$, is said to be a
{\sl polar  oriented closed} manifold.

\noindent
ii) If the  volume form $\mu$ on a compact complex manifold $M$
is nonvanishing and meromorphic with only first order poles
on a normal crossing divisor $N\subset M$, then $M$ is called
a {\sl polar oriented} manifold {\sl with boundary}.
The hypersurface $N$ is then endowed with a polar
orientation $\nu:=\res\,\mu\neq0$ and $(N,\nu)$ is called
the {\sl polar boundary} of $(M,\mu)$.
\end{defn}

\begin{rem}
By definition, polar oriented closed manifolds are complex manifolds whose
canonical bundle is trivial (Calabi--Yau, Abelian manifolds or, for example,
any complex tori, if we do not restrict ourselves to algebraic manifolds).
We have just defined the notion of the polar orientation in a more
restrictive sense than before, when we considered the definition of
chains. In fact, polar chains with their orientations are to be compared
to oriented piece-wise smooth submanifolds in differential topology,
while the ambient space on which we want to have Poincar\'e duality has
to be everywhere smooth and oriented.
Zeros of a volume form could be regarded as a complex analogue
of singularities  of a real manifold.\footnote
{For instance, on a complex curve $X$ of genus $g$
one has $HP_1(X)=\C^{\,g}$,
and a holomorphic 1-differential representing a generic element in
$HP_1(X)$ has $2g-2$ zeros. From this point of view the complex
genus $g$ curve is like a
graph which has $g$ loops joined by $g-1$ edges and having
$2g-2$ trivalent (i.e., ``non-smooth'') points. The ``smooth cases'' are
$\C\P^1$, which corresponds to a real segment, and an elliptic curve, which is
a complex counterpart of the circle in this precise sense.}

\end{rem}

\medskip
{\sl a) Polar intersection number}

\smallskip

\begin{ssect}{}\label{Xnumber}
Let $(M,\mu)$ be a polar oriented  closed manifold of dimension $n$.
In such a case we define the following  natural pairing between
its polar homology groups $HP_p(M)$ and $HP_{n-p}(M)$
of complimentary dimension.

According to Proposition~\ref{hom},  the above groups can
be mapped to the Dolbeault cohomology groups
$H^{n,n-p}_{\bar\partial}(M)$ and $H^{n,p}_{\bar\partial}(M)$, respectively.
On a polar oriented closed manifold we are given a nowhere vanishing
section $\mu$ of the line bundle $\Omega^n_M\:$.
Hence, we have the  isomorphism
$H^{n,n-p}_{\bar\partial}(M)\stackrel{\mu^{-1}}{\longrightarrow}
H^{0,n-p}_{\bar\partial}(M)$.
Using this and the product in Dolbeault cohomology we obtain the
following pairing:
$$
H^{n,p}_{\bar\partial}(M) \otimes H^{n,n-p}_{\bar\partial}(M)
\xrightarrow{{\rm id}\otimes\mu^{-1}}
H^{n,p}_{\bar\partial}(M) \otimes H^{0,n-p}_{\bar\partial}(M)
\to H^{n,n}_{\bar\partial}(M)
\stackrel{\!\!\sim}{\to} \C \;.
$$
Together with the homomorphism
$\rho\!: HP_{k}(X)\to H^{n,n-k}_{\bar\partial}(X)$,
this yields the pairing

$$
\langle~\Cdot~\rangle\!: HP_p(M) \otimes HP_{n-p}(M)
\to \C \;.
$$
Here, in fact, we interchanged the order of factors (see
the explicit formula (\ref{Xintegral}) below).

\end{ssect}

\begin{ssect}{}\label{Xintegral}
Consider two polar cycles, $[a]\in HP_p(M)$ and $[b]\in HP_{n-p}(M)$.
Let $t_a$ and $t_b$ be the Dolbeault forms representing
$\rho([a])$ and $\rho([b])$ respectively. Then, the above product, which we
denoted by $\langle a\Cdot b\rangle$, can be written as
$$
\langle a\Cdot b\rangle = \int_M\: t_b \,\wedge\,\frac{t_a}{\mu} \;.
$$
Note that $t_a$ is an $(n,n-p)$-form and thus, $t_a/\mu$
is a $(0,n-p)$-form that can be integrated against an $(n,p)$-form
$t_b\,$.

\end{ssect}

\begin{defn} The pairing $\langle a\Cdot b\rangle$ of polar cycles is
called the {\sl polar intersection index}.
\end{defn}

\begin{rem}
If the conjecture (\ref{Conj}) is true,
this pairing is non-degenerate.

\end{rem}

\begin{ssect}{}
Let us consider now the case when the cycles $a$ and $b$ are smooth
and transverse. That is $a=(A,\alpha)$ and $b=(B,\beta)$, where $A$ is
a smooth $p$-dimensional variety and $\alpha$ a holomorphic $p$-form
on it (and similarly for $(B,\beta)$ in dimension $n-p$) and it is
assumed that $A$ and $B$ intersect transversely. Then, we have
the following formula for the polar intersection index.

\end{ssect}

\begin{th}\label{Xtrans}
The polar intersection index of two smooth transverse cycles
$(A,\alpha)$ and $(B,\beta)$ is given by the following sum over the
set of points in $A\cap B$:
$$
\langle(A,\alpha)\Cdot (B,\beta)\rangle = \sum_{P\in A\cap B}
\frac{\alpha(P)\wedge\beta(P)}{\mu(P)}
$$
Here $\alpha(P)$ and $\beta(P)$ are understood
as exterior forms on $T_P M = T_P A \times T_P B$ obtained by the
pull-back
from the corresponding factors.
\smallskip

{\rm
The ratio in the right-hand-side can be understood as the comparison
of the polar orientations brought to the intersection point $P$ by the
two cycles with the polar orientation $\mu(P)$ of the ambient
manifold at that point.}

\end{th}

\begin{proof}{}
As we have already mentioned,
the homomorphism $\rho$ of Sect.~\ref{hom} can be
conveniently described in terms of the following natural map of polar chains:
$$
\varphi:\, {\cal C}_k(M) \to {\cal D}^{n,n-k}(M) \;,
$$
where ${\cal D}^{p,q}(M)$ is the space of currents of degree $(p,q)$
(i.e., a space of  linear
functionals on smooth $(n-p,n-q)$-forms, see
\cite{GH}). As a matter of fact, this map is
described by the integral (\ref{pair}).
For a $p$-dimensional submanifold $A\subset M$,
let the current $\delta_A\in{\cal D}^{n-p,n-p}(M)$  denote  the
linear functional on $(p,p)$-forms corresponding to
the integration over $A$.
The current $\delta_A$ is supported on $A$.
Therefore, for a $p$-form $\alpha$ defined on $A$, the product
$\delta_A\wedge\alpha$
makes sense and defines a current in ${\cal D}^{n,n-p}(M)$.
Recalling the
isomorphism of the cohomology of currents with the cohomology of
smooth forms,
$$
H^j({\cal D}^{i,\bullet}(M),\bar\partial)
\stackrel{\sim\,}{\longrightarrow}
H^{i,j}_{\bar\partial}(M) \;,
$$
we can use $\delta_A\wedge\alpha$ and $\delta_B\wedge\beta$ in place of
$t_a$ and $t_b$ in the integral (\ref{Xintegral}). Thus, for
a transverse intersection of smooth polar cycles we derive
that\footnote
{Although the (exterior) product of currents is not in general defined, the
product of the cohomology classes of $\bar\partial$-closed currents is always
well
defined. In some cases it is easy to find a representative of the product of
cohomology classes. For instance, for two submanifolds $V$ and $W$ in a
generic position the product of cohomology classes of the corresponding
currents $\delta_V$ and $\delta_W$ is represented by $\delta_{V\cap W}$. In
this sense we can write $\delta_V\wedge\delta_W=\delta_{V\cap W}$.}
$$
\langle(A,\alpha)\Cdot (B,\beta)\rangle =
\int_M\: \delta_B\wedge\beta \,\wedge\, \frac{\delta_A\wedge\alpha}{\mu}
= \int_M\: \left(\frac{\alpha\wedge\beta}{\mu}\right)
\cdot \delta_A\wedge\delta_B \;.
$$
The second equality can be checked in local coordinates.
This proves the theorem, since $\delta_A\wedge\delta_B$ is supported
on $A\cap B$.
\end{proof}

\medskip
{\sl b) Polar intersection product}
\smallskip

\begin{ssect}{}\label{Lambda}
Now consider the case when on a polar oriented  closed manifold $(M,\mu)$
we have two polar cycles of arbitrary dimensions $p$ and $q$
(not necessarily complimentary ones).
Similarly to the pairing (\ref{Xnumber}),
we may consider the following chain of homomorphisms:
\begin{eqnarray*}
HP_p(M) \otimes HP_{q}(M)
\xrightarrow{s_{12}\circ(\rho\otimes\rho)}
H^{n,n-q}_{\bar\partial}(M) \otimes H^{n,n-p}_{\bar\partial}(M) \\
\xrightarrow{{\rm id}\otimes\mu^{-1}}
H^{n,n-q}_{\bar\partial}(M) \otimes H^{0,n-p}_{\bar\partial}(M)
\to H^{n,2n-p-q}_{\bar\partial}(M) \;,
\end{eqnarray*}
where $s_{12}$ is the transposition of tensor factors and
the last term is understood as zero unless $p+q\geqslant n$.
Let $\Lambda$
denote the resulting composition:
$$
\Lambda\!: HP_p(M) \otimes HP_{q}(M) \to
H^{n,2n-p-q}_{\bar\partial}(M) \;.
$$

\end{ssect}

\begin{ssect}{}
If the conjecture (\ref{Conj}) holds,
the above homomorphism and the
inverse of $\rho$ in (\ref{hom}) define an intersection product on
polar homology,
$$
HP_p(M) \otimes HP_{q}(M) \to HP_{p+q-n}(M) \;.
$$
However, even without this hypothesis
we will show that if $a$ and $b$ are two smooth
transverse cycles, $[a]\in HP_p(M)$, $[b]\in HP_{q}(M)$,
their polar intersection product can be represented in $HP_{p+q-n}(M)$
by a smooth cycle $c$.

\end{ssect}

\begin{ssect}{$\C\,$-orientations of vector spaces.}
Let $W$ be an $n$-dimensional complex vector space, and
$\mu$ be a non-zero complex volume form on $W$ ($\mu\in\bigwedge^nW^*$). Let
$V_A,V_B\subset W$ be vector subspaces of dimensions
$\dim V_A=p,\,\dim V_B=q,\,p+q\geqslant n$, which intersect
transversely, i.e.\ $V_A+V_B=W$ (or, $r:=\dim V_A\cap V_B=p+q-n\,$).
Suppose we are also given
complex volume forms on each of $V_A$ and $V_B$, that is
$\alpha\in\bigwedge^pV_A^*$ and $\beta\in\bigwedge^qV_B^*$.
(We may say that
all three spaces $W, V_A$, and $V_B$ are $\C\,$-oriented.) Then, one can
define a complex volume form $\gamma$ on  the intersection $V_A\cap V_B$
(i.e., one can $\C\,$-orient it) as follows.

Let $\lambda_A\in\bigwedge^{n-p}W^*$ be a non-zero exterior form
conormal to $V_A$ (i.e., $\lambda_A$ vanishes on any vector from $V_A$
and is non-zero as an element of $\bigwedge^{n-p}(W/V_A)^*$). Similarly, let
$\lambda_B\in\bigwedge^{n-q}W^*$ be a non-zero exterior
form conormal to $V_B$. Note that in this case $\lambda_A\wedge\lambda_B$
is a form conormal to $V_A\cap V_B$.

\end{ssect}

\begin{deflem} \label{lam1}
Given complex orientations (i.e., complex
volume forms) $\alpha, \beta$, and $\mu$ of the vector subspaces $V_A,V_B$
and the space $W$ respectively, the following complex orientation $\gamma$ of
the intersection $V_A\cap V_B$ is defined by the following relation:
$$
          \lambda_A\wedge\lambda_B\wedge\gamma =
          \left(\frac{\lambda_A\wedge\alpha}{\mu}\right)\cdot
          \left(\frac{\lambda_B\wedge\beta}{\mu}\right)
          \cdot\mu  \;.
$$
Here $\alpha,\beta$, and $\gamma$ are understood as arbitrary
extensions of these forms to the whole space $W$. The $r$-form $\gamma$
on $V_A\cap V_B$ depends neither on these extensions, nor on the choice
of the auxiliary forms $\lambda_A$ and $\lambda_B$.
\end{deflem}

\begin{proof}{} A straightforward verification.
\end{proof}

\begin{corll} \label{lam2}
For a transverse intersection of subspaces of complimentary dimensions
($p+q=n$ and  $W=V_A\oplus V_B$), the 0-form $\gamma$ is just the following
complex number:
$$
               \gamma=\frac{\alpha\wedge\beta}{\mu} \;,
$$
where $\alpha$ and $\beta$ are now understood as exterior forms on
$W=V_A\times V_B$ obtained by pull-backs from the corresponding factors.
\end{corll}

\begin{ssect}{}\label{c=ab}
Let $a=(A,\alpha)$ and $b=(B,\beta)$ be two smooth polar cycles of
dimension $p$ and $q$ respectively in a polar oriented closed manifold
$(M,\mu)$, $p+q\geqslant n=\dim M$. Suppose they intersect
transversely. Then we can define a $(p+q-n)$-cycle $c=(C,\gamma)$,
where $C=A\cap B$ is a smooth subvariety in $M$ and $\gamma$ is a
holomorphic $(p+q-n)$-form on it defined by (\ref{lam1}).
Let us denote this as $[c]=[a]\Cdot[b]$ and call the intersection
product.
(Of course, the product $[a]\Cdot[b]$ equals zero if $p+q<n$.)
\end{ssect}

\begin{th}\label{Lambda=rho}
The polar intersection product $[a]\Cdot[b]$ of two smooth transverse
cycles $a=(A,\alpha)$ and $b=(B,\beta)$ defined above agrees with the
homomorphism {\rm(\ref{Lambda})}, i.e.,
$$
\Lambda([a]\otimes[b]) = \rho([a]\Cdot[b]) \;,
$$
where $\rho: HP_{k}(X)\to H^{n,n-k}_{\bar\partial}(X)$
was defined in {\rm(\ref{pair})-(\ref{hom})}.
\end{th}

\begin{proof}{}
This will be similar to the proof of Theorem~\ref{Xtrans} and will use the
same notations. We first represent the polar cycles $a$ and $b$ by the
currents $\delta_A\wedge\alpha$ and $\delta_B\wedge\beta$
respectively, then
$$
\Lambda([a]\otimes[b]) = \left[\delta_B\wedge\beta\wedge\:
\frac{\delta_A\wedge\alpha}{\mu} \right] \;,
$$
where $[~]$ on the right is understood as taking the
$\bar\partial$-cohomology class. On the other hand, for $c=(C,\gamma)$
introduced in (\ref{c=ab}), the current representing $c$ is
$\delta_C\wedge\gamma$ and it is easy to show that
$$
\delta_C\wedge\gamma =
\delta_B\wedge\beta\wedge\:
\frac{\delta_A\wedge\alpha}{\mu} \;,
$$
what implies the statement of the Theorem. The last equality is easily
checked by noticing that $\delta_A$ is an $(n-p,n-p)$-form (in fact, a
current) conormal to $A$ and similarly for $\delta_B$, while
$\delta_C=\delta_A\wedge\delta_B\,$. This is to be compared to
$\lambda_A$ and $\lambda_B$ in (\ref{lam1}).
One has to note only that, e.g., $\delta_A$ is conormal to $A$
over $\R$ (that is in the sense of $(n-p,n-p)$-forms) while
$\lambda_A$ is conormal to it over $\C$ (that is in the
sense of $(n-p,0)$-forms).
\end{proof}

\begin{rem}\label{X-sgranicej}
We have defined the polar intersection on any complex manifold $M$ that
can be equipped with a {\it holomorphic} non-vanishing volume form $\mu$.
This is analogous to the topological intersection theory on a compact smooth
oriented manifold without boundary. (Note that the Poincar\'e duality in
this context should correspond to the Serre duality.)
Furthermore, the consideration above easily extends to the case of
a complex manifold
possessing a {\it meromorphic} non-vanishing form $\mu$
(in particular, to a complex projective space),
i.e., to the case of a polar oriented manifold $(M,\mu)$ with boundary
$(N,\res\,\mu)$ (cf. \ref{boundary-def}).
The latter setting is similar
to the topological intersection theory on manifolds with boundary.
In this case the above formulas can be used to define the pairing between
polar homology $HP_k(M)$ and polar homology relative to the boundary
$HP_{n-k}(M,N)$.

\end{rem}



\bigskip

\begin{ackn}
B.K. and A.R. are grateful for hospitality to the Max-Planck-Institut f\"ur
Mathematik in Bonn where this work was conceived and later completed.
B.K. is also grateful to the Institute for Advanced Study in Princeton,
and A.R. is also grateful to the Erwin Schr\"odinger International
Institute for Mathematical Physics in Vienna
where a part of this work was done.
The present work was partially sponsored by PREA and McLean awards.

A.R. is extremely grateful to A.~Bondal and
A.~Levin for explaining some valuable tools of algebraic geometry
and  very
instructive discussions. A number of valuable remarks by A.~Levin
saved us from many missteps in the course of this work.

The work of B.K. was partially supported by an Alfred P. Sloan
Research Fellowship,   by the NSF and NSERC research grants.
The work of A.R. was supported in part by the Grants RFBR-98-01-00344,
INTAS-97-0103 and the Grant 00-15-96557 for the support of scientific
schools.

\end{ackn}


\begin{smallbibl}{MMMM}

\bibitem[AKSZ]{Sch}\ M.~Alexandrov, M.~Kontsevich, A.~Schwarz,
and O.~Zaboronsky,\\
{\em The Geometry of the Master Equation and
Topological Quantum Field Theory}\/,\\
Int. J. Mod. Phys. {\bf A12} (1997) 1405--1430
~(hep-th/9502010)

\bibitem[A]{Arn}\ V.I.~Arnold,
{\em  Arrangement of ovals of real plane algebraic curves,
 involutions
of smooth four-dimensional  manifolds, and on arithmetic
of integral-valued quadratic forms}\/,
Func. Anal. and Appl.  {\bf 5:3}  (1971),  169--176

\bibitem[BO]{BO}\ S.~Bloch and A.~Ogus,
{\em Gersten's conjecture and the
homology of schemes}\/, \\
Ann. Sci. \'Ecole Norm. Sup. (4) {\bf 7} (1974), 181--201.

\bibitem[De]{De}\ P.~Deligne, {\em Th\'eorie de Hodge. II}\/, Inst. Hautes
\'Etudes Sci. Publ. Math. {\bf 40} (1971), 5--57.

\bibitem[DT]{DT}\ S.K.~Donaldson and R.P.~Thomas,
{\em Gauge theory in higher
dimensions}\/. \\ The geometric universe (Oxford, 1996),
Oxford Univ. Press, Oxford (1998) 31--47.

\bibitem[FK]{FK}\ I.B.~Frenkel and B.A.~Khesin,
{\em Four Dimensional Realization of Two Dimensional Current
Groups}\/, Commun. Math. Phys. {\bf 178} (1996) 541--561

\bibitem[FT]{FT}\ I.B.~Frenkel and A.N.~Todorov,
paper in preparation

\bibitem[Ger]{Ger}\ A.~Gerasimov, unpublished (1995)

\bibitem[Gr]{Gr}\ P.A.~Griffiths,~ {\em Variations on a Theorem of
Abel}\/, Invent. Math., {\bf 35} (1976) 321--390

\bibitem[GH]{GH}\ P.A.~Griffiths and J.~Harris, {\em Principles of
Algebraic Geometry}\/, Wiley, NY (1978)

\bibitem[KR]{KR}\ B.~Khesin and A.~Rosly, {\em Symplectic geometry on
moduli spaces of holomorphic bundles over complex surfaces}\/,
In: The Arnoldfest, Proceedings of a conference in honour
of V.I.~Arnold for his sixtieth birthday, -- Toronto 1997,
Editors: E.~Bierstone et al.
Fields Institute/AMS Communications (1999) {\bf 24}, 311--323
[math.AG/0009013]

\bibitem[iKR]{iKR}\ B.~Khesin and A.~Rosly,~ in preparation

\bibitem[ASL]{ASL}\ A.S.~Losev, private communication (1999)

\bibitem[LNS]{arhar}\ A.~Losev, N.~Nekrasov, and S.~Shatashvili,
{\em Issues in Topological Gauge Theory}\/, \\
Nucl. Phys. {\bf B534} (1998) 549--611
~(hep-th/9711108)

\bibitem[T]{T}\ R.P.~Thomas,
{\em Gauge theory on Calabi--Yau manifolds}\/, Ph.D. thesis, Oxford
(1997), 1--104

\bibitem[Tju]{Tju} Tjurina, G. N.
{\em On the moduli space of complex surfaces with $q=0$ and $K=0$}\/,
Chapter IX in {\em Algebraic surfaces}\, by Shafarevich, I. R. et al.,
Trudy Steklov Math. Inst. {\bf 75} (1965) 163--191

\bibitem[W]{W}\ E.~Witten,
{\em Chern--Simons gauge theory as a string theory}\/,
The Floer memorial volume, Progr. Math., {\bf 133}, Birkh\"auser,
Basel (1995) 637--678 ~(hep-th/9207094)

\end{smallbibl}

\end{document}